\documentclass[letterpaper, 10 pt, conference]{IEEEconf}  
\IEEEoverridecommandlockouts                              
\usepackage{multirow, array}
\usepackage{epsfig}
\usepackage{algorithmicx}
\usepackage{algpseudocode}
\usepackage{cite}
\usepackage{amsmath,amssymb,amsfonts,mathptmx,colortbl,braket}
 
\usepackage[outdir=./]{epstopdf}
\usepackage{amsthm,BOONDOX-cal}
\usepackage{graphicx}
\usepackage{enumerate}
\usepackage{epsfig}
\usepackage[font=small,labelfont=bf]{caption}
\usepackage{amsmath} 
\usepackage{epsfig}
\usepackage{epstopdf}
\usepackage{algorithm,array}
\usepackage{algpseudocode}
\usepackage{multirow,multicol}
\let\labelindent\relax
\usepackage{subfig,braket,enumitem,adjustbox,booktabs,slashbox}
\usepackage{color}
\usepackage{float}
\newfloat{algorithm}{t}{lop}
\usepackage{pgf,tikz}
\usetikzlibrary{arrows}
\usepackage{cuted}
\usepackage[overload]{empheq}
 \usepackage{tabularx, booktabs, makecell, caption}
    \usepackage{siunitx}

\newcommand{\oprocendsymbol}{\hbox{$\bullet$}}
\newcommand{\oprocend}{\relax\ifmmode\else\unskip\hfill\fi\oprocendsymbol}

\newtheorem{theorem}{Theorem}[section]

\newtheorem{lemma}[theorem]{Lemma}

\newtheorem{remark}[theorem]{Remark}

\renewcommand{\theenumi}{(\roman{enumi})}

\newcommand{\real}{\ensuremath{\mathbb{R}}}

\newcommand{\Hc}{\mathcal{H}}

\newcommand{\Kc}{\mathcal{K}}
\newcommand{\Lc}{\mathcal{L}}

\newcommand{\Oc}{\mathcal{O}}

\newcommand{\Zc}{\mathcal{Z}}

\newcommand{\T}{^{\top}}

\newcommand{\zeros}{\mathbf{0}}

\newcommand\tr{\operatorname{tr}}

\newcommand\pd{\partial}

\newcommand\td{\text{d}}

\newcommand\Kbe{\Kc_{\beta}}
\newcommand\grad{\nabla}

\newcommand\gradK{\nabla_K}

\newcommand\Pin{P_{\infty}}

\newcommand{\Hci}{\mathcal{H}_{\infty}}
\newcommand\st{\text{s.t.}}

\newcommand{\longthmtitle}[1]{\mbox{}{\textit{(#1).}}}
\newcommand\Hinf{\Hc_{\infty}}
\newcommand\fe{f_{\eta}}

\DeclareMathOperator*{\argmax}{arg\,max}

\usepackage{enumitem}

\newlist{notes}{enumerate}{1}
\setlist[notes]{label=Note: ,leftmargin=*}

\renewcommand{\baselinestretch}{0.97}
\allowdisplaybreaks

\begin{document}
\title{A Scalable Procedure for $\Hci-$Control Design}
\author{Amit Kumar$^{*}$ and Prasad Vilas Chanekar
\thanks{$^{*}$Corresponding author.} 
\thanks{This research is supported by SERB Start-up Research Grant - SRG/2023/001636.}
\thanks{AK and PVC are with the Department of Electronics and
		Communication Engineering, Indraprastha Institute of Information Technology, New Delhi, India 110020, {\tt \{amitku,\,prasad\}}@iiitd.ac.in}}
	\maketitle

    \begin{abstract}
        This paper proposes a novel gradient based scalable procedure for $\Hci-$control design. We compute the gradient using algebraic Riccati equation and then couple  it  with a novel Armijo rule inspired step-size selection procedure. We perform numerical experiments of the proposed solution procedure on an exhaustive list of benchmark engineering systems to show its convergence properties. Finally we compare our proposed solution procedure with available semi-definite programming based  gradient-descent algorithm to demonstrate its scalability.
    \end{abstract}

\section{Introduction}

A practical engineering system is influenced by various modelling uncertainties and external disturbances, which can significantly compromise its overall performance, stability and reliability.
Robust control is a specialized area within control theory that focuses on designing controllers capable of preserving performance amid disturbances and uncertainties.
The effectiveness of a robust controller is typically evaluated using a $\gamma$ parameter, which quantifies the impact of model uncertainties and external disturbances on the system performance. A smaller $\gamma$ value signifies a more reliable and robust system behaviour \cite{skogestad2007multivariable}.

 $\Hci$-control is a specialized form of robust control that addresses the impact of disturbances on system output in a worst-case scenario~\cite{zhou1996robust,stoorvogel1992h}. 
 %
The $\Hci$-control problem in the input-output framework was formally introduced in \cite{zames1981feedback}. Early methods for solving the $\Hci-$control problem include interpolation methods \cite{limebeer1988interpolation}, frequency domain methods \cite{francis1987course}, polynomial methods \cite{kwakernaak1986polynomial}, spectral factorization methods \cite{kimura1989conjugation}, and time-domain optimal control methods \cite{tadmor1990worst}, among others.
In \cite{doyleglover1989}, a groundbreaking approach was taken to develop an observer-based sub-optimal $\Hci$ controller by effectively solving two algebraic Riccati equations (AREs) \cite{zhou1996robust}. 
In \cite{khargonekar1991h_} did a thorough investigation into the $\Hci$ control problem, addressing critical challenges such as input disturbances, measurement noise, and uncertainties in the initial system state for both finite and infinite time horizons. 
A novel one-shot design process was introduced to derive a sub-optimal closed-form solution for $\Hci$-control, which is parameterized by the variable $\gamma$. This innovative approach is specifically tailored for a class of systems that meet certain geometric criteria, as outlined in \cite{saberi1994closed}. In the context of stable systems, the parameter $\gamma$ plays a crucial role by effectively establishing the upper limit of the $\Hci$-norm \cite{zhou1996robust,stoorvogel1992h}. An efficient bisection algorithm to compute the parameter $\gamma$  was presented in \cite{boyd1989bisection}.

Within the  state space framework, the $\Hci$ control design problem was approached as a norm minimization problem that includes a stability constraint \cite{apkarian2017h}. Typically, these problems involve a scalar objective function in conjunction with nonlinear matrix equality or inequality constraints, resulting in a computationally challenging non-convex formulation \cite{balakrishnan2003semidefinite, boyd1994linear, apkarian2017h}. The non-convex, nonlinear constraint is expressed as a bilinear matrix inequality (BMI), which can be transformed into a convex linear matrix inequality (LMI)~\cite{balakrishnan2003semidefinite,boyd1994linear} with the help of a suitable transformation. The convex transformed problem was solved using standard convex optimization techniques \cite{boyd2004convex} resulting a globally optimal solution.
However the magnitude of entries in the resultant optimal gain was large entries \cite{scherer1989h} making it infeasible for practical use. In \cite{gahinet1994linear} a more useful strategy involved the development of sub-optimal controllers for $\Hci$-control using an Linear Matrix Inequalities (LMI) based convex optimization approach  for both continuous and discrete-time systems.
Necessary and sufficient conditions for the existence of $\Hci$ controllers, articulated through the framework of LMIs, were introduced in \cite{iwasaki1994all}. In addition, an  iterative method which employed a sequential semi-definite programming approach with local super-linear convergence properties to address the LMI-constrained $\Hci$ control design problem was developed in \cite{fares2002robust}. 
Recently \cite{ma2022symmetric} rigorously investigated the  $\Hci$-control design problem by utilizing the iterative alternating direction method of multipliers (ADMM) \cite{boyd2011distributed} in conjunction with the sophisticated symmetric Gauss-Seidel technique \cite{adegbege2016gauss}. In \cite{ICC} a novel gradient-based formulation which uses LMIs was proposed for $\Hci$-control design.

The existing literature indicates that convex optimization based design methods are often associated with high gain values. Conversely, classical approaches, such as frequency domain analysis and spectral factorization, frequently prove inadequate for large systems. LMI-based iterative methods suffer for scalability issues. With respect to the aforementioned research gaps we  outline the contributions of our paper. 

\textit{Paper Contribution:} We first present refinements in the semi-definite programming (SDP) based problem formulation for the $\Hci$-control design problem in our recent work \cite{ICC}. We analyze and point out its shortcomings with regard to time complexity. We then propose the $\Hci$-control optimization problem with an algebraic Riccati equation constraint. We then propose a novel iterative gradient-based solution procedure that uses a novel gradient based on the upper bound of the original problem. We also propose a new Armijo rule inspired \cite{Bertsekas} step-size selection procedure with a dynamic parameter. Note that the controller gain matrix synthesized using our algorithm has the same order of magnitude as the entries of the initial gain matrix. We perform numerical analysis of our proposed solution procedure with regards to performance and scaling properties on an exhaustive list of benchmark practical engineering examples \cite{compleib_models}. Finally, we show the efficacy of our proposed algorithm by comparing its performance and scalability with the existing gradient-based solution procedure that uses SDPs. 

The paper is organized as follows: Section II presents the Preliminaries followed by the Problem Formulation in  Section III. Sections IV and Section V presents Solution Procedure and Examples followed by concluding remarks in Section VI.

\textit{Notations:} We let $\real$  denote the set of real numbers. For a matrix $X\in \real^{p\times q}$ having $p-$rows and $q-$columns, we use $X^{\top}, \ \tr(X),\ \Vert X\Vert _{\infty},\ \Vert X\Vert\ \text{and}\ \sigma_{\max}(X) $ for its transpose, its trace, its $\Hci-$ norm, its Frobenius norm and its largest singular value respectively. For a vector $x\in \real^n$, $\Vert x\Vert_2$ denotes its Euclidean norm. We use $I$ and $\zeros$ to denote the identity matrix and matrix of zeros, respectively, of appropriate sizes. For $X\in \real^{a\times a}$, $X(\succeq)\succ0$ denotes $X$ is positive (semi-)definite. For scalar function $f$,  notation $\grad_X f$ denotes the gradient of $f$ with respect to $X$. For a time-dependent signal $\phi(t)$, $\Vert \phi(t)\Vert_2^2 = \int_0^{\infty}\phi\T(t)\phi(t)\, \td t.$
    
\section{Preliminaries}\label{prelim}
    Consider a linear time-invariant (LTI) dynamical system 
\begin{align} \label{dynamics1}
    \notag \dot{x} &= Ax + Bu + B_1w, \quad x(0) = x_0,\\
    z&=Cx+Du,\quad u = Kx,
 \end{align}
where $A \in \mathbb{R}^{n\times n}$, $B \in \mathbb{R}^{n\times m}$, $B_1 \in \mathbb{R}^{n\times m_1}$, $C \in \mathbb{R}^{r\times n}$ and $D \in \mathbb{R}^{r\times m}$ are the system matrices. Whereas $x \in \mathbb{R}^{n}$, $u \in \mathbb{R}^{m}$, $w \in \mathbb{R}^{m_1}$ and $z \in \mathbb{R}^{r}$ denote the state, the full-state feedback control input, the exogenous disturbance input and the controlled output of the system respectively. $K \in \mathbb{R}^{m\times n}$ denotes the full-state feedback controller gain matrix. We make the following assumptions for the system \eqref{dynamics1}. 
\begin{enumerate}
    \item $(A,\:B)$ and $(A,\:B_1)$ are stabilizable.
    \item $(A,\:C)$ is detectable.
\end{enumerate}
 Let \( G(s) \) be the closed-loop transfer function  between the input disturbance and the output. Then we have,
\begin{align*}
    G(s) = C_1(sI-A_c)^{-1}B_1,
\end{align*}
where $A_c=A+BK,\;C_1=C+DK$. Using $G$, the worst-case performance metric is defined using the $\Hci$-norm of the transfer function  in  time and frequency domains as follows~\cite{zhou1996robust}.
\begin{align} \label{definition of H-infinity}
    \Vert G\Vert _{\infty}= \sup_{w\neq 0} \frac{\Vert z\Vert_2}{\Vert w\Vert_2}= \sup_{\omega \in \mathbb{R}} \sigma_{\text{max}}(G(j\omega)).
\end{align}
 The transfer function \( G \) depends on the controller gain \( K \). The \(\Hci\)-control problem for the system  \eqref{dynamics1} involves finding a stabilizing gain \( K \) that minimizes the worst-case performance objective \( \Vert G\Vert_{\infty} \) and is represented by the following optimization problem:
\begin{align}
    \label{problem1}
    \notag &\min_{K,\; \gamma >0} \quad \gamma\\
   & \st \qquad \Vert G(s) \Vert_{\infty} < \gamma,\\
     & \notag \quad \qquad A+BK \;\; \text{is Hurwitz.}
\end{align}
Next, the problem \eqref{problem1} is suitably formulated in a computationally viable format.

\section{Problem Formulation}
\label{problem-formulation}
 From semi-definite programming theory \cite{balakrishnan2003semidefinite,boyd1994linear}, \eqref{problem1} can be represented as a non-convex matrix optimization problem with bi-linear matrix inequality (BMI) constraint as follows,
\begin{align}
    \label{nonconvex-problem1}
    \notag &\min_{K,\; \beta >0,\; P\succ 0} \quad \beta\\
   & \st \qquad  \begin{pmatrix}
        A_c\T P + PA_c +C_1\T C_1 & PB_1\\
        B_1\T P & -\beta I
    \end{pmatrix}  \preceq \zeros.
\end{align}
Note that here $\gamma = \sqrt{\beta}$. The objective function in \eqref{nonconvex-problem1} is linear, so is convex. The non-convexity arises from the BMI constraint due to the product between $P$ and $K$. 

For a known stabilizing $K$ (i.e., $A_c=A+BK$ Hurwitz) the problem \eqref{nonconvex-problem1} is formulated as
\begin{align}
    \label{primal1}
    \notag &\min_{\beta_K >0,\; P\succ 0} \quad \beta_K\\
   & \st \qquad  \begin{pmatrix}
        A_c\T P + PA_c +C_1\T C_1 & PB_1\\
        B_1\T P & -\beta_K I
    \end{pmatrix}  \preceq \zeros.
\end{align}
\eqref{primal1} is a convex optimization problem with $A_c$ known and is known as the `\textit{primal problem}'. 

The dual of  problem \eqref{primal1} is as follows \cite{balakrishnan2003semidefinite},
\begin{align}
      \label{dual1}
    \max_{Z} \qquad &
    \notag J_K = \tr (C_1Z_{11}C_1\T) \\
    \text{s.t.} \qquad & Z_{11}A_c\T + A_cZ_{11} + Z_{12}B_1\T + B_1Z_{12}\T = \zeros,\\
    \notag & \begin{pmatrix}
        Z_{11} & Z_{12} \\
        Z_{12}\T & Z_{22}
    \end{pmatrix} \succeq  \zeros ,\quad \tr(Z_{22}) = 1.
    \end{align}
As $K$ is known, the \textit{`dual problem'} \eqref{dual1} is a concave problem with a matrix equality constraint.

 Let $\beta^*_K$ and $J^*_K$ be the optimal values of \eqref{primal1} and \eqref{dual1} respectively and let $\gamma_K$ be the $\Hinf-$norm of the system for the gain $K$. 
  From duality theory \cite{boyd2004convex,balakrishnan2003semidefinite}  we have $\beta^*_K\geq J^*_K.$  However, from our numerical experiments performed on the benchmark practical test problems in control \cite{compleib_models}, we observe strong duality for any random stabilizing gain $K$ (i.e., $A_c$ Hurwitz)  \cite{boyd2004convex} i.e.,
  $$\beta^*_K=J^*_K=\gamma_K^2.$$ A rigorous proof of this observation will be provided in the journal extension of this paper. Once a stabilizing $K$ is known, $\gamma_K$ is computed using the bisection method as described in \cite{boyd1989bisection}. We now define the following sets. 
  \begin{align}
      \label{set-def1}
      &\Kc = \Set{K|K\in \real^{m \times n},\; A_c=A+BK \text{  is Hurwitz}},\\
     \notag& \Zc_K=\Set{Z| \begin{array}{ll}
     K\in \Kc, \;A_c=A+BK\\ Z=\begin{pmatrix}
        Z_{11} & Z_{12} \\
        Z_{12}\T & Z_{22}
    \end{pmatrix} \succeq  \zeros,\, \tr(Z_{22})=1,\\
      Z_{11}A_c\T + A_cZ_{11} + Z_{12}B_1\T + B_1Z_{12}\T = \zeros
\end{array}}, \\
\notag &\Zc^* = \Set{Z^* |K \in \Kc,\; Z^* = \argmax_{Z \in \Zc_K} \;\tr (C_1Z_{11}C_1\T) }.
  \end{align}
Obviously $\Zc^*\subset\Zc_K.$ Using the strong duality observation and \eqref{set-def1} we can now rewrite \eqref{nonconvex-problem1} as follows
   \begin{align}
     \label{prob1}
        \min_{K\in \Kc,\; Z^* \in \Zc^*} \quad J^*_K.
 \end{align}
  Here $J^*_K=\tr (C_1Z_{11}^*C_1\T)$. Note that for each $K\in\Kc$ a unique $Z^*$ exists in set $\Zc_K$  which is computed by solving the convex optimization problem \eqref{dual1}. The set $\Zc^*$ is the collection of all such $Z^*\in \Zc_K$ for each $K\in \Kc$. The gradient of $J^*_K$ with respect to $K$ at $K=K^0$ is given by \cite{ICC},
  \begin{align} \label{gradient expression}
    \grad_K J^*_K(K^0) = 2(B\T L + D\T C + D\T D K^0)Z^*_{11}, 
\end{align} 
with $A_c= A+BK^0$ and
 \begin{subequations}
    \begin{align*}
    Z^*_{11}A_c\T + A_cZ^*_{11} + Z^*_{12}B\T + B {Z^*}_{12}\T &= \zeros, \\
     LA_c+ A_c\T L + C_1\T C_1+M_{11} &=\zeros,\\  
    \tr(Z^*_{22})=1,\quad LB_1+M_{12} =\zeros,\quad \nu I+M_{22} &=\zeros, \\     Z^* \succeq \zeros,\; L=L\T\in\real^{n\times n},\; MZ^*=\zeros, \; M &\succeq \zeros.
\end{align*}
\end{subequations}
 Once $K^0$ is known $Z^*,L,M,\nu$ are computed using \eqref{dual1} and semi-definite programming (SDP) solvers like $CVX$ \cite{cvx}. Using \eqref{gradient expression}, we can now solve \eqref{prob1}  using gradient-descent algorithm coupled with Armijo step-size selection procedure \cite{Bertsekas}. However, in the gradient-descent procedure applied to \eqref{prob1}, one needs to solve SDP \eqref{dual1} multiple times. The time-complexity of solving an SDP using solvers for a $n-$state system is $\Oc(n^6)$ \cite{gramlich2023structure}. Thus as the system state size $n$ increases, it becomes computationally infeasible to solve the $\Hci-$norm controller synthesis problem.

Next, we propose a novel problem formulation and scalable solution procedure for $\Hci-$norm controller synthesis.
\section{Solution Procedure}
\begin{figure}[h]
    \centering
    \includegraphics[width=0.5\linewidth, height=0.3\linewidth]{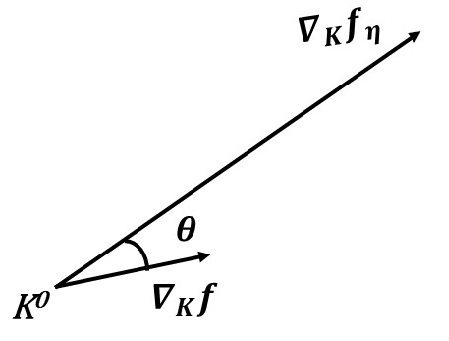}
    \caption{$\gradK f$ and $\gradK \fe$ at $K^0\in\Kbe$}
    \label{fig1}
\end{figure}

For a known $K$ and $A_c$ Hurwitz,  $\gamma_K^2=\beta_K$  is knowm from the bisection procedure \cite{boyd1989bisection}.  Then  the algebraic Riccati equation (ARE) $$A_c\T \Pin +\Pin A_c + C_1\T C_1 + \frac{1}{\beta}\Pin B_1 B_1^T\Pin= \zeros,$$ has a solution $\Pin\succeq 0$ for  all $\beta\geq \beta_K$  \cite[Chapter 16]{zhou1996robust}. The problem \eqref{nonconvex-problem1} can be rewritten as
\begin{align}
    \label{prob2}
    \notag &\min_{K,\; \beta >0,\; \Pin\succ 0} \quad f=\beta\\
   & \st \qquad  
        A_c\T \Pin +\Pin A_c + C_1\T C_1 + \frac{1}{\beta}\Pin B_1 B_1^T\Pin= \zeros.
\end{align}
Problem \eqref{prob2} is nonconvex and difficult to solve as it is even when $K$ is known due to the quadratic ARE constraint. Consider the set definition
\begin{align*}
& \Kbe=\Set{K| \begin{array}{ll}
     K\in \Kc, \; \exists \,\beta >0,\; \exists\Pin\succ \zeros,\\
      A_c\T \Pin +\Pin A_c + C_1\T C_1 + \frac{1}{\beta}\Pin B_1 B_1^T\Pin= \zeros.
\end{array}}.
  \end{align*}
  \eqref{prob2} can be rewritten as 
     \begin{align}
     \label{prob4}
        \min_{K\in \Kbe} \quad f(K)=\beta.
 \end{align}
 Now we can apply gradient-descent algorithm to iteratively solve \eqref{prob4}. However we try to compute the $\gradK f$ we counter the following condition in the process
$$\Big(A+BK^0+\frac{1}{\beta}B_1B_1\T\Pin\Big)L_0 + L_0\Big(A+BK^0+\frac{1}{\beta}B_1B_1\T \Pin\Big)\T  = \zeros.$$ Above has only the trivial solution $\Pin= \zeros.$  This makes application of gradient-based methods to solve \eqref{prob2} impossible. To circumvent the aforementioned shortcoming we propose the following problem
\begin{align}
    \label{prob3}
    \notag &\min_{K,\; \beta >0,\; \Pin\succ 0} \quad \fe =\beta+\eta \tr(\Pin)\\
   & \st \qquad  
        A_c\T \Pin +\Pin A_c + C_1\T C_1 + \frac{1}{\beta}\Pin B_1 B_1^T\Pin= \zeros,
\end{align}
with $\eta >0$. Now $$f\leq\fe.$$ However for a sufficiently small value of $\eta>0$, $\fe$ is a tight upper-bound of $f$. In such situation we conjecture that $f$ and $\fe$ behave almost the same as $K\in \Kbe$ varies. Hence to solve \eqref{prob4} using gradient-descent algorithm we use  $\gradK \fe$ as the search direction at $K^0\in\Kbe$ instead of $\gradK f$. Note that in computation of $\gradK \fe$ the condition of $\Pin = \zeros$ is avoided as shown in the next result. To derive our result, we will make use of known matrix properties~\cite{lewis2012optimal}. For matrices $X,\: Y$,  $\tr(X\:Y) = \tr(Y\:X)$, $\tr(X) = \tr(X\T)$ provided the matrices are compatible for multiplication, and $\frac{\pd\:\tr(Y\:X) }{\pd\:X}= Y\T$.
\begin{lemma}\longthmtitle{Gradient of $\fe$ with respect to $K$}\label{lemma-gradient}
Consider optimization problem  \eqref{prob3}. Then for $K^0\in \Kbe$, the gradient of $\Lc$ with respect to $K$  is
\begin{align}\label{grad2}
    \gradK \fe(K^0)=2(B\T P_\infty + D\T C + D\T D K^0)L_0, 
          \end{align}
with
 \begin{subequations}\label{eq-mas}
    \begin{align}
    A_c\T P_\infty +P_\infty A_c + C_1\T C_1 + \frac{1}{\beta}P_\infty B_1 B_1^TP_\infty= \zeros,\label{eq-mas1}\\
     A_1L_0 + L_0A_1\T  +\eta I= \zeros,\label{eq-mas2}
          \end{align}
\end{subequations}
where $A_1 = A_c+\frac{1}{\beta}B_1B_1\T\Pin$.
\end{lemma}
\begin{proof}
Consider the Lagrangian function of \eqref{prob3}, 
\begin{align*}
 \Lc(K,\Pin,L_0) &= \beta+\eta\tr(\Pin)+\\
      &\hspace{0cm}\tr\big(L_0(A_c\T P_\infty +P_\infty A_c + C_1\T C_1 + \frac{1}{\beta}P_\infty B_1 B_1^TP_\infty)\big).
\end{align*}
$\Lc$ is a function of $K,\Pin,L_0$. From the definition of $\Kbe$, for every $K\in\Kbe$ we have $\fe(K) =\Lc(K)$. Thus $\fe$ and $\Lc$ have the same behavior as $K$ varies in $\Kbe$. Therefore at $K^0\in \Kbe$ we have $\gradK \fe(K^0) = \gradK \Lc(K^0).$ Now $\frac{\pd \Lc}{\pd L_0}$ gives \eqref{eq-mas1},  $\frac{\pd \Lc}{\pd \Pin}$ gives \eqref{eq-mas2} and
 $\frac{\pd \Lc}{\pd K}$ gives \eqref{grad2}.
 %
\end{proof}
From our numerical experiments we observe that at any $K^0\in\Kbe$, $\gradK f$ and $\gradK \fe$ are aligned at an acute angle i.e., $0< \theta \ll \frac{\pi}{2}$ as shown in Fig \ref{fig1}. However $\Vert \gradK \fe \Vert \gg \Vert \gradK f \Vert.$
Due to this the direct use of  $\gradK \fe$ in place of $\gradK f$ in the gradient-descent algorithm creates issues in finding step-size using Armijo rule \cite{Bertsekas}. Hence we  make suitable modifications in the Armijo rule while computing the step-size. We next present our gradient-based solution algorithm. The input to the algorithm is $K^0\in\Kbe$ and $f(K^0) = \beta^0$ computed using bisection method \cite{boyd1989bisection}.
\begin{algorithm}[H]
\caption{Gradient-descent algorithm} \label{alg:Newton}
\begin{algorithmic}[1]
\Require $A,\;B,\;\eta,\; K^{0},\;\beta^0,\;\gradK \fe(K^0),\;\epsilon$ 
\Ensure $K^*,\; \beta^*$
\State Set $j\gets0$, $K^j\gets K^0$
 \While{1}
\State Compute step size $s^j$ using \textit{modified-Armijo rule}\Comment{See Algorithm~\ref{alg:Armijo}}
\State $K^{j+1}\gets K^j-s^j\gradK \fe(K^j))$
\State Compute $f(K^{j+1}) =\beta^{j+1}$ using bisection method \cite{boyd1989bisection}
\If {$\vert f(K^{j+1})-f(K^j)\vert\leq \epsilon$}
                \State Go to Step 12
                \EndIf
\State $j\gets j+1$,\; $K^j\gets K^{j+1}, \;\beta^j=\beta^{j+1}$
\State Compute $\gradK \fe(K^j)$ using Lemma \ref{lemma-gradient}
\EndWhile
\State $K^*\gets K^j, \; \beta^*\gets \beta^j$ 
\State \Return $K^*,\;\beta^*$
\end{algorithmic}
\end{algorithm}
The original version of Armijo rule has fixed parameters \cite{Bertsekas}. 
However we make one of the parameters dynamic and state the modified-Armijo rule to compute the step-size as follows.
\begin{algorithm}[h]
		\caption{modified-Armijo rule}\label{alg:Armijo}
		\begin{algorithmic}[1]
			\Require$ A,\; B,\; \eta,\;  \alpha^0,\; \zeta,\;  K^{j},\; f(K^j), \;\gradK \fe(K^j) $ 
			\Ensure $s^j$
            \State Set $\alpha \gets \alpha^0$
             \State $m^j = \gradK \fe(K^j)$
                \State Set $s^j\gets 1$
                \While{$f(K^j- s^j m^j)\geq f(K^j) -\alpha s^j\big\Vert m^j \big\Vert^2 $}
                \State{~$s^j\gets\zeta s^j$,}
                \If {$s^j < 10^{-15}$}
                \State $s^j \gets 1$, $\alpha \gets \alpha/5$
                \EndIf
                \State Compute $f(K^j- s^j m^j)$ using bisection method \cite{boyd1989bisection}
                \EndWhile
   \State \Return $s^j$
			\end{algorithmic}	
	\end{algorithm} 
    
Typically in Algorithm \ref{alg:Armijo}, $\alpha^0 = 0.3$ and $\zeta = 0.5$. 
Note that in Algorithm \ref{alg:Armijo} to obtain $f(K^j- s^j m^j)$ we again use the bisection algorithm \cite{boyd1989bisection}  at $K=K^j- s^j m^j$. In Algorithm \ref{alg:Newton}, the modified-Armijo method ensures that the objective function $i.e., f(K^j) = \beta^j$ keeps on decreasing as the iteration count $j$ increases. As the objective function is lower bounded by zero hence we must have convergence to a stationary point \cite{dennis1996numerical}. Also, for a sufficiently small $\epsilon$, the stopping criterion (Step 12) $\vert f(K^{j+1})-f(K^j)\vert\leq \epsilon$,  may also ensure $\gradK f(K^j) \approx 0$. From our numerical experiments we observe this behavior for any initial $K^0\in \Kbe$ hence we conjecture existence of a global convergence to stationary point phenomenon. 
 \begin{remark}\longthmtitle{Time-complexity comparison of gradient-descent methods using SDP and ARE}
		{\rm 
		As stated at the end of Section \ref{prelim}}, the gradient-descent method using SDP has a time-complexity of $\Oc(n^6)$ per iteration. While our proposed Algorithm \ref{alg:Newton} solves one ARE and one Lyapunov equations to compute the gradient per iteration. The time-complexity of this process is $\Oc(n^3)$ per iteration \cite{gramlich2023structure}. Also the bisection method to compute $f=\beta$ is inherently fast. Thus our proposed algorithm scales very well   in comparison with the algorithms who use SDPs. This in also verified by numerical experiments presented in the next section.  
		\oprocend
\end{remark}
 A detailed analytical study regarding the convergence and scaling properties of Algorithm \ref{alg:Newton} along with the modified-Armijo rule is part of our future work. Next, we present practical numerical examples to justify our proposed controller synthesis procedure.

\section{Examples} \label{examples}
 We  demonstrate the utility of our proposed algorithm \ref{alg:Newton} through two examples. Example 1 is reproduced from \cite{9844674} while Example 2 is a set of benchmark real-world examples from the COMP$l_e$ib library \cite{compleib_models}. We compare the performance between the gradient-descent algorithm with SDP formulation (\eqref{prob1} with gradient \eqref{gradient expression}) and our proposed gradient-descent algorithm, which uses AREs (Algorithm 1). We use the improvement in the initial and final $\gamma$ as the performance metric along with time durations of the corresponding solution procedures. The following parameters are used to implement Algorithm \ref{alg:Newton}: $\epsilon =  10^{-5}, \ \alpha^0 = 0.3, \ \zeta = 0.5,\ \eta = 0.1$. The initial input $K^0$ is the linear-quadratic regulator (LQR controller)  \cite{lewis2012optimal} computed with the weighing matrices $Q = I$ and $R = I$ of appropriate sizes.
 
In these examples, the algorithm is implemented in MATLAB R2024a and executed on a laptop with an Intel(R) Core(TM) i7-12700H processor (14 Cores), 2.30 GHz speed and 16 GB of RAM. We use MATLAB \cite{MATLAB} for all simulations and programming purposes. To solve the convex problem \eqref{dual1} in Algorithm \ref{alg:Newton}, we use CVX solver \cite{cvx} with `SDPT3' semi-definite programming package \cite{toh1999sdpt3,tutuncu2003solving}.
\subsection{Example $1$}
Consider the dynamics in  \cite[Example 2]{9844674} where the model has parametric uncertainties. The nominal values for system matrices are as follows.
\begin{align*}
    &{A} = \begin{pmatrix}
        0.2229 & 0.5637\\
    0.8708 & 0.9984 \\
    \end{pmatrix}, \; {B}  = \begin{pmatrix}
        0.5254 & 0.6644\\
        0.3872 & 0.9145\\
    \end{pmatrix}\\
    &{B_1} = \begin{pmatrix}
        1 &0 \\
    0 &1 \\
    \end{pmatrix} ,\;  {C} = \begin{pmatrix}
        1 & 0 \\
    0 & 1 \\
    0 & 0 \\
    0 & 0 \\
    \end{pmatrix}, \; {D} = \begin{pmatrix}
        0 & 0 \\
        0 & 0 \\
        1 & 0 \\
        0 & 1 \\
    \end{pmatrix}.
\end{align*}
The optimized controller design in \cite{9844674} is $K = -\begin{pmatrix}
        0.9643 & 2.1060 \\
        0.2088 & 5.6843 \\
    \end{pmatrix}$ and with the corresponding $\gamma = 4.9411$. However Algorithm \ref{alg:Newton} in this paper computes $K^* = -\begin{pmatrix}
        0.8426 & 0.9893  \\
        0.0551 & 2.5743  \\
    \end{pmatrix}$ with the corresponding $\gamma^* = 2.6736$  almost half the $\gamma$ determined in \cite[Example 2]{9844674} showing the efficacy of our approach.   Note that the entries of $K^0$ and $K^*$ have same magnitude order.

\subsection{Example $2$}
We apply our proposed Algorithm \ref{alg:Newton} to the set of engineering benchmark examples from the COMP$l_e$ib library \cite{compleib_models}. This collection encompasses systems of various sizes. We present a comparative analysis of our proposed Algorithm 1 (ARE) \ref{alg:Newton}, alongside gradient-descent procedure using SDPs \cite{ICC} in the accompanying Table \ref{table}. Each system in Table \ref{table} is accompanied with the state size $n$ and number of inputs $m$. As both procedures start at the same initial $K^0$, the optimal $\gamma^*$ from our proposed procedure matches or betters the optimal $\gamma^*$ from the SDP-based procedure. However, the computation time of our proposed Algorithm \ref{alg:Newton} is significantly less compared to the SDP-based procedure. For large-size systems with $n=60,\; 120,\; 240,\; 480,\; 960$, the SDP-based procedure failed to produce any solution while our proposed Algorithm \ref{alg:Newton} designed a controller with improved performance. This shows the efficacy of our proposed $\Hci-$control design procedure. Note that the  norms for initial and final gain matrices for all the benchmark examples i.e., $\Vert K^0\Vert$ and $\Vert K^*\Vert$ have the same magnitude order. 
\begin{table*}[t] 
    \centering
    \renewcommand{\arraystretch}{1.3} 
    \begin{tabular}{|c|c|c|c|c|c|c|c|c|c|}
        \hline
        \textbf{System} & \textbf{$n$} & \textbf{$m$} & \textbf{$\gamma^0$-initial} & \textbf{$\gamma^*$-SDP} & \textbf{$\gamma^*$-ARE} & \multicolumn{2}{c|}{\textbf{Time (seconds)}} & \textbf{\% improvement in $\gamma$-SDP} & \textbf{\% improvement in $\gamma$-ARE} \\
        \cline{7-8}
        & & & & & & \textbf{SDP} & \textbf{ARE} & & \\
        \hline
        AC1 & 5 & 3 & 0.05 & 0.03 & 0.03 & 23.22 & 1.10 & 40.00 & 40.00 \\
        AC2 & 5 & 3 & 0.18 & 0.11 & 0.11 & 18.71 & 1.26 & 38.89 & 38.89 \\
        AC3 & 5 & 2 & 3.87 & 3.68 & 3.70 & 43.54 & 1.15 & 4.91 & 4.39 \\
        AC4 & 4 & 1 & 6.43 & 1.07 & 1.11 & 16.28 & 1.06 & 83.36 & 82.74 \\
        AC8 & 9 & 1 & 2.11 & 1.58 & 1.57 & 288.03 & 2.65 & 25.12 & 25.59 \\
        AC11$^*$ & 5 & 2 & 3.94 & 2.92 & 2.92 & 54.95 & 1.30 & 25.89 & 25.89 \\
        AC12 & 4 & 3 & 3.63 & 2.54 & 2.13 & 68.94 & 4.78 & 30.03 & 41.32 \\
        AC17 & 4 & 1 & 6.79 & 6.61 & 6.61 & 6.01 & 0.76 & 2.65 & 2.65 \\
        AC18 & 10 & 2 & 32.93 & 4.66 & 4.49 & 74.44 & 7.02 & 85.86 & 86.37 \\
        HE1$^*$ & 4 & 2 & 0.20 & 0.07 & 0.06 & 30.29 & 1.03 & 65.00 & 70.00 \\
        HE2$^*$ & 4 & 2 & 3.86 & 2.60 & 2.59 & 49.39 & 0.95 & 32.64 & 32.90 \\
        HE3 & 8 & 4 & 0.99 & 0.84 & 0.85 & 11.07 & 1.87 & 15.15 & 14.14 \\
        HE4 & 8 & 4 & 14.11 & 12.99 & 12.98 & 69.93 & 1.67 & 7.93 & 8.04 \\
        HE5 & 8 & 4 & 2.62 & 2.18 & 2.17 & 16.49 & 1.37 & 16.79 & 17.18 \\
        REA1$^*$ & 4 & 2 & 1.06 & 0.68 & 0.65 & 29.43 & 3.59 & 35.85 & 38.68 \\
        REA2$^*$ & 4 & 2 & 1.07 & 0.68 & 0.63 & 26.27 & 2.30 & 36.45 & 41.12 \\
        DIS1$^*$ & 8 & 4 & 5.36 & 4.28 & 4.28 & 88.58 & 2.44 & 20.15 & 20.15 \\
        DIS2 & 3 & 2 & 1.22 & 0.95 & 0.93 & 649.94 & 6.11 & 22.13 & 23.77 \\
        DIS4 & 6 & 4 & 1.58 & 1.14 & 0.96 & 183.81 & 17.51 & 27.85 & 39.24 \\
        DIS5 & 4 & 2 & 54.98 & 44.50 & 44.40 & 86.09 & 1.53 & 19.06 & 19.23 \\
        AGS & 12 & 2 & 8.17 & 8.17 & 8.17 & 9.78 & 1.58 & 0.00 & 0.00 \\
        BDT1 & 11 & 3 & 0.29 & 0.27 & 0.27 & 11.65 & 3.35 & 6.90 & 6.90 \\
        MFP$^*$ & 4 & 3 & 6.61 & 4.28 & 4.15 & 62.45 & 1.40 & 35.26 & 37.21 \\
        IH & 21 & 11 & 12.71 & 4.34 & 3.58 & 529.43 & 38.14 & 65.85 & 71.83 \\
        EB1 & 10 & 1 & 2.01 & 1.90 & 1.90 & 168.74 & 2.92 & 5.47 & 5.47 \\
        EB2$^*$ & 10 & 1 & 0.76 & 0.51 & 0.50 & 212.65 & 3.14 & 32.89 & 34.21 \\
        EB3$^*$ & 10 & 1 & 0.76 & 0.53 & 0.50 & 239.66 & 4.19 & 30.26 & 34.21 \\
        TF2$^*$ & 7 & 2 & 1.12 & 0.29 & 0.28 & 57.21 & 2.83 & 74.11 & 75.00 \\
        TF3$^*$ & 7 & 2 & 1.12 & 0.29 & 0.28 & 62.31 & 2.80 & 74.11 & 75.00 \\
        PSM & 7 & 2 & 0.93 & 0.92 & 0.92 & 3.40 & 1.06 & 1.08 & 1.08 \\
        NN1 & 3 & 1 & 19.45 & 13.18 & 13.18 & 821.87 & 1.03 & 32.22 & 32.22 \\
        NN2 & 2 & 1 & 1.92 & 1.52 & 1.52 & 553.17 & 3.05 & 20.83 & 20.83 \\
        NN4 & 4 & 2 & 1.94 & 1.38 & 1.38 & 41.54 & 1.27 & 28.87 & 28.87 \\
        CM1$^*$ & 20 & 1 & 1.00 & 0.97 & 0.90 & 114.17 & 19.30 & 3.00 & 10.00 \\
        TMD$^*$ & 6 & 2 & 5.02 & 2.61 & 2.57 & 81.16 & 2.19 & 48.01 & 48.80 \\
        CM2$^{**}$ & 60 & 1 & 1.00 & $\times$ & 0.88 & $\times$ & 268 & $\times$ & 12 \\
        CM3$^{**}$ & 120 & 1 & 1.00 & $\times$ & 0.90 & $\times$ & 805 & $\times$ & 10 \\
        CM4$^{**}$ & 240 & 1 & 1.00 & $\times$ & 0.90 & $\times$ & 4766 & $\times$ & 10 \\
        CM5$^{**}$ & 480 & 1 & 1.00 & $\times$ & 0.90 & $\times$ & 9383 & $\times$ & 10 \\
        CM6$^{**}$ & 960 & 1 & 1.00 & $\times$ & 0.91 & $\times$ & 46420 & $\times$ & 9 \\
        \hline
    \end{tabular}
    \caption{Analysis of COMP$l_e$ib benchmark examples. $\% \; \text{improvement} = \frac{\gamma^0-\gamma^*}{\gamma^0}\;\times \;100$. $^*$ systems have $\alpha^0 = 0.5 \ \& \ \epsilon = 10^{-9}. $ $^{**}$ systems have $\alpha^0 = 0.3 \ \& \ \epsilon = 10^{-15}.$ $\times $ means no solution obtained.}
    \label{table}
\end{table*}

\section{Conclusion}
 This paper presents an iterative gradient-based procedure for the $\Hci-$ control design problem which used a novel gradient expression based on the upper bound of the original problem. A modified Armijo rule inspired step-size selection rule is also proposed. The convergence and scalability properties of the proposed control design algorithm were numerically analyzed by applying it to an exhaustive list of benchmark engineering systems. A comparison with the convex optimization-based approach available in the literature was made to show the efficacy of the proposed theory. Future work involves analytical analysis of the convergence and scalability properties of the proposed solution procedure.

\bibliographystyle{IEEEtran}
\bibliography{ref1}
    
\end{document}